\documentclass[amsmath,twocolumn,secnumarabic,superscriptaddress,%
    floatfix,amssymb,aps]{revtex4}
\usepackage{times,graphicx,amsthm}

\graphicspath{{./figs/}}
\newcommand{\intersect}{\cap}
\newcommand{\union}{\cup}
\newcommand{\R}{\mathbb{R}}
\renewcommand{\H}{\mathbb{H}}
\renewcommand{\d}{\partial}
\newcommand{\clip}{\text{cl}}
\renewcommand{\sb}{\text{sb}}
\newcommand{\setm}{\smallsetminus}
\newcommand{\h}{h} 
\def\nthhull/{{$n^{\text{th}}$ hull}}

\newtheorem*{maintheorem}{Main Theorem}
\newtheorem*{theorem}{Theorem}
\newtheorem{lemma}{Lemma}
\newtheorem*{corollary}{Corollary}
\newtheorem*{proposition}{Proposition}
\newtheorem*{conjecture}{Conjecture}
\theoremstyle{definition}
\newtheorem*{definition}{Definition}

\begin{document}

\title{The Second Hull of a Knotted Curve}
\date{May 5, 2000; revised February 1, 2003} 

\author{Jason Cantarella}
\email[Email: ]{cantarel@math.uga.edu}
\affiliation{Department of Mathematics, University of Georgia,
Athens, GA 30602}

\author{Greg Kuperberg}
\email[Email: ]{greg@math.ucdavis.edu}
\affiliation{Department of Mathematics, University of California,
Davis, CA 95616}

\author{Robert B. Kusner}
\email[Email: ]{kusner@math.umass.edu}
\affiliation{Department of Mathematics, University of Massachusetts,
Amherst, MA 01003}

\author{John M. Sullivan}
\email[Email: ]{jms@math.uiuc.edu}
\affiliation{Department of Mathematics, University of Illinois,
Urbana, IL 61801}

\def\incgrh#1#2{\includegraphics[height=#2]{figs/#1}}

\begin{abstract}
The convex hull of a set $K$ in space consists of points which are,
in a certain sense, ``surrounded'' by $K$.
When $K$ is a closed curve, we define
its higher hulls, consisting of points which are ``multiply surrounded''
by the curve.  Our main theorem shows that if a curve is knotted then
it has a nonempty second hull.  This provides a new proof of the
F{\'a}ry/Milnor theorem that every knotted curve has total curvature at
least~$4\pi$.
\end{abstract}

\maketitle

A space curve must loop around at least twice to become knotted.
This intuitive idea was captured in the celebrated F{\'a}ry/Milnor theorem,
which says the total curvature of a knotted curve $K$ is at least $4\pi$.
We prove a stronger result:
that there is a certain region of space doubly enclosed,
in a precise sense, by $K$.  We call this region the second hull of $K$.

The convex hull of a connected set $K$ is characterized by the fact
that every plane through a point in the hull must interesect $K$.
If $K$ is a closed curve, then a generic plane intersects $K$
an even number of times, so the convex hull is
the set of points through which every plane cuts $K$ twice.
Extending this, we introduce a definition which measures the degree to
which points are enclosed by~$K$: a point $p$ is said to be in the
\emph{\nthhull/} of $K$ if every plane through $p$ intersects $K$ at
least $2n$ times. 

A space curve is convex and planar if and only if it has the
\emph{two-piece property}: every plane cuts it into at most two pieces. This
has been generalized by Kuiper (among others; see \cite{Kuiper}) in the
study of tight and taut submanifolds. It was Milnor
\cite{milnor} who first observed that for a knotted curve, there are
planes in every direction which cut it four times.

\begin{figure}[ht]
  \begin{center}
    \parbox{2.0in}{\rotatebox{22}{\incgrh{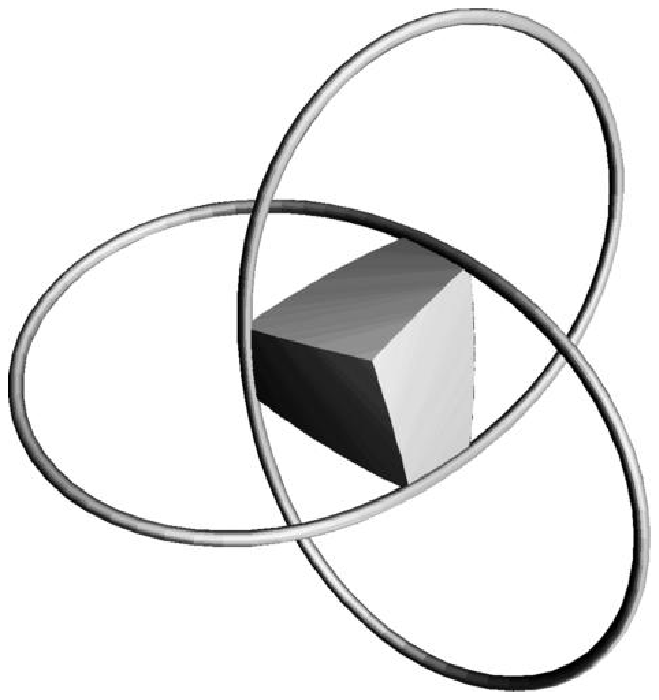}{115bp}}} 
    \hspace{0.0in}
    \parbox{1.25in}{\vspace{-0.25in} \rotatebox{270}{\includegraphics[height=44bp]{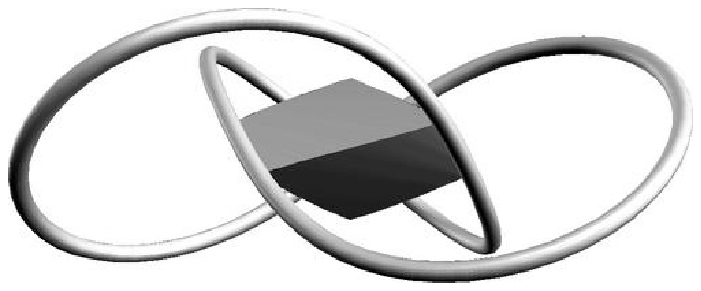}}} 
  \end{center}
  \vspace{-0.5in}
  \caption[The second hull of a trefoil.]
    {This figure shows the second hull of a trefoil knot, the shaded
    region in the center of the knot, from two perspectives.
    Every plane which intersects second hull must cut the knot at least
    four times.
    The trefoil shown here is a polygonal approximation to the critical
    point---and presumed minimizer---for the
    M\"obius-invariant knot energy of Freedman, He and Wang~\cite{FHW}.
    The region shown is an accurate computation of the second
    hull of this polygonal knot, produced by John Foreman.  The sharp
    edges remind us that second hulls, like convex hulls, need not
    be as smooth as the curves which generate them.
} \label{fig:trefoil2hull}
\end{figure}

Our main theorem strengthens this to find points through which every
plane cuts the knot four times.  More precisely, it says
that any curve $K$ with a knotted inscribed polygon (in particular
any knotted curved of finite total curvature)
has a nonempty second hull,
as shown for a trefoil in Figure~\ref{fig:trefoil2hull}.
Furthermore, in a certain sense, this hull must extend across the knot:
Using quadrisecants, we get an alternate approach to second hulls,
and in particular we show that the second hull of a
generic prime knot $K$ is not contained in any embedded normal tube around~$K$.

We use some ideas from integral geometry to show that the total curvature
of a space curve with nonempty \nthhull/ is at least $2\pi\:\!n$.
Thus, our theorem provides an alternate proof of the
F\'ary/Milnor theorem.

We do not yet know any topological knot type which is guaranteed
to have a nonempty third hull, though we suspect that
this is true for the (3,4)--torus knot (see Figure~\ref{fig:8_19}).
It would be interesting to find, for $n\ge3$, computable topological
invariants of knots which imply the existence of nonempty \nthhull/s.

Our definition of higher hulls is appropriate only for curves,
that is, for one-dimensional submanifolds of $\R^d$.  Although we can
imagine similar definitions for higher-dimensional submanifolds,
we have had no reason yet to investigate these.  For $0$-dimensional
submanifolds, that is, finite point sets $X\subset\R^d$,
Cole, Sharir and Yap~\cite{CSY} introduced the notion of
the \emph{$k$-hull} of~$X$, which is the
result of removing from $\R^d$ all half-spaces containing
fewer than $k$ points of $X$.  Thus it is the set of points
$p$ for which every closed half-space with $p$ in the boundary
contains at least $k$ points of $X$.
Although we did not know of their work when we first formulated
our definition, our higher hulls for curves are similar in spirit.

\section{Definitions and Conventions}

\begin{definition}
Let $K$ be a closed curve in $\R^d$,
by which we mean a continuous map from the circle, modulo reparametrization.
The \emph{\nthhull/} $\h_n(K)$ of $K$
is the set of points $p\in\R^d$ such that $K$ cuts every
hyperplane through $p$ at least $2n$ times, where nontransverse
intersections are counted according to the conventions below.
\end{definition}

Here, we basically count the intersections of $K$ with a plane~$P$ by counting
the number of components of the pullback of $K \cap P$ to the domain circle.
If the intersections are transverse, then (thinking of~$P$
as horizontal, and orienting $K$) we find equal numbers of upward and downward
intersections; the total number of intersections (if nonzero)
equals the number of components of $K \setm P$ (again counted in the domain).

To handle nontransverse intersections, we again orient $K$
and we adopt the following
conventions.  First, if $K\subset P$, we say $K$ cuts~$P$ twice (once in
either direction). If $K\cap P$ has infinitely many components, then we
say $K$ cuts~$P$ infinitely often.  Otherwise, each connected
component of the intersection is preceded and followed by open arcs in
$K$, with each lying to one side of~$P$.  An \emph{upward} intersection
will mean a component of $K\cap P$ preceded by an arc below~$P$
\emph{or} followed by an arc above~$P$.  (Similarly, 
a \emph{downward} intersection will mean a component preceded
by an arc above~$P$ or followed by an arc below~$P$.)
A glancing intersection, preceded and followed by arcs on
the same side of~$P$, thus counts twice, as both an upward and a
downward intersection.

With these conventions, it is easy to see that
an arc whose endpoints are not in a
plane~$P$ intersects $P$ an even number of times if the endpoints are
on the same side of~$P$ and an odd number of times otherwise.  For a
closed curve, there are always equal numbers of intersections in each
direction. 

We will later apply the same definition to disconnected curves $L$, like links.
In that context, of course, the number of intersections of $L$ with $P$
cannot simply be counted as the number of components of $L\setm P$.

\section{The Second Hull of a Knotted Curve}

In order to prove that any knot has a second hull,
we start with some basic lemmas about \nthhull/s.

\begin{lemma}\label{lem:cvx}
For any curve~$K$, each connected component
of its \nthhull/ $\h_n(K)$ is convex.
Furthermore, $\h_n(K)\subset\h_{n-1}(K)$.
\end{lemma}
\begin{proof}
Let $C$ be a connected component of $\h_n(K)$.  If $q$ is a point in
its convex hull, then by connectedness, any plane~$P$ through $q$ must
meet $C$ at some point $p$.  Thus $P$ is also a plane through $p\in
C\subset \h_n(K)$, so $P$ meets $K$ at least $2n$ times.
By definition, then, $q\in\h_n(K)$.  The second
part of the statement then follows from the definition of~$\h_n$.
\end{proof}

A useful lemma allows us to simplify a given curve, with control on
its higher hulls.

\begin{lemma}\label{lem:cut}
Let $K^*$ be the result of replacing some subarc $\Gamma$ of $K$ by the
straight segment $S$ connecting its endpoints $p$, $q$. No plane cuts
$K^*$ more often than it cuts $K$, so $\h_n(K^*)\subset \h_n(K)$.
\end{lemma}
\begin{proof}
Let $P$ be any plane.
If $K^*\cap P$ has infinitely many components then so does $K\cap P$.  So
we may assume $K^*\cap P$ has only a finite number of components.
The generic case is when $p,q\notin P$. If they are on the same side
of $P$, then $S$ misses $P$ and there is nothing to prove. Otherwise,
$S$ has a unique intersection with $P$, but $\Gamma$ must also have cut $P$
at least once.

The other cases require us to invoke our conventions for counting
intersections, and thus require a bit more care.  First suppose
$p,q \in P$. If $K^* \subset P$, this counts as two intersections,
but by our conventions $K$ cuts $P$ a nonzero even number of times,
hence at least twice. Otherwise, let $A$ be the component of $K^* \cap P$
containing $S$.  It counts as one or two intersections with $P$,
depending on whether points just before and after $A$ are on
different sides of $P$ or the same side.
But the corresponding arc of~$K$ then cuts
$P$ a nonzero odd or even number of times respectively, hence at
least as often as $A$ does.

Next, suppose $p\in P$ and $q$ is above $P$.  Let $A$ be the
component of $K\cap P$ containing $p$, and consider the arc of~$K$
from a point $r$ just before $A$ to $q$.  This arc cuts $P$ a nonzero
odd or even number of times if $r$ is below or above $P$, but then the
corresponding arc in $K^*$ intersects $P$ once or twice, respectively.

In any case, $K^*$ has at most as many intersections as $K$
did with the arbitrary plane $P$.
The final statement of the lemma then follows from the definition of $\h_n$.
\end{proof}

We will need to locate regions of space
which can be said to contain the knottedness of $K$.

\begin{definition}
We say that an (open) half-space $H$ \emph{essentially contains} a
knot $K$ if it contains $K$, except for possibly a single
point or a single unknotted arc.  
\end{definition}

Note that if $H$ essentially contains $K$, then $K$ cuts $\d H$ at
most twice.  Thus, we can \emph{clip $K$ to $H$}, replacing the single
arc (if any) outside $H$ by the segment in $\d H$ connecting its
endpoints.  The resulting knot $\clip_H(K)$ is isotopic to~$K$,
and is contained in $\overline H$.

Our key lemma now shows that essential half-spaces remain essential
under clipping:

\begin{lemma}\label{lem:clip}
Let $K^*:=\clip_H(K)$ be the clipping of a knot $K$ to some
half-space $H$ which essentially contains it.
Then any half-space $H'$ which essentially
contains $K$ will also essentially contain $K^*$.
\end{lemma}
\begin{proof}
We must show that either $H'$ contains $K^*$, or the single arc of this
new knot outside $H'$ is unknotted.
If $H'$ contains $K$ then, 
since $K^*$ lies within the convex hull of $K$,
$H'$ must contain $K^*$, and we are done.

\begin{figure*}[ht]
  \begin{center}     \parbox{2in}{\includegraphics{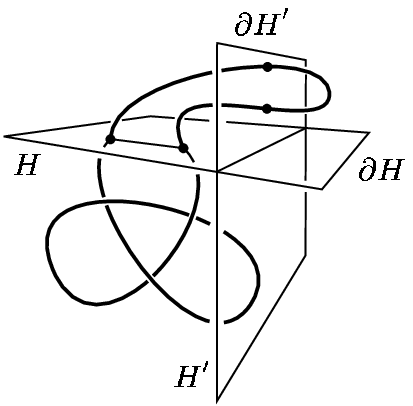}}
    \hspace{0.125in} \parbox{2in}{\includegraphics{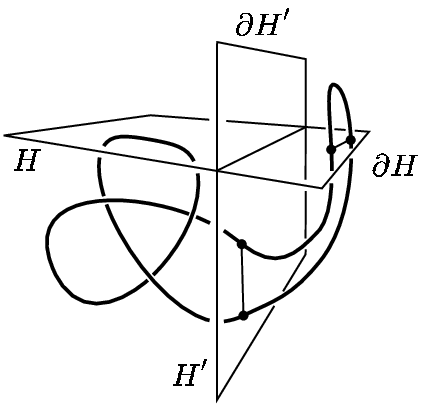}}
    \hspace{0.125in} \parbox{2in}{\includegraphics{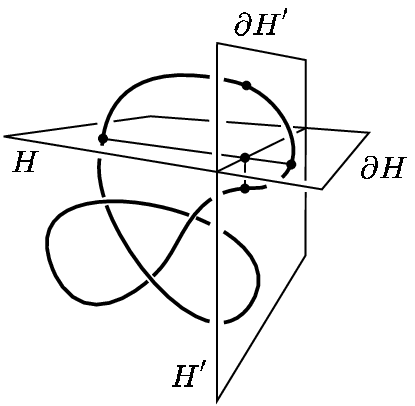}}
  \end{center} 
\caption[The cases of Lemma~\ref{lem:clip}.]
    {This figure shows the three cases of interest in the proof of
    Lemma~\ref{lem:clip}.  On the left, we see the case where both
    intersections of $K$ with $\d H'$ occur outside $H$. In the
    center, we see the case where both intersections of $K$ with
    $\d H'$ occur \emph{inside} $H$. On the right, we see the case
    where~$K$ intersects~$\d H'$ once inside and once outside~$H$.
    In each figure, the segments in~$\d H$ and~$\d H'$ replacing
    unknotted subarcs of $K$ when $K$ is clipped to these planes
    are drawn in lightly. The half-space $H$ is below the
    horizontal plane $\d H$, while the half-space $H'$ is to the
    left of the vertical plane $\d H'$.
\label{fig:clip}}
\end{figure*} 

Otherwise, consider the single arc of $K$ outside $H'$, and count how
many of its endpoints are outside $H$.
The different cases are illustrated in Figure~\ref{fig:clip}.
If both endpoints are outside~$H$, again we find $H'$ contains $K^*$.
If neither is, the unknotted arc of~$K$ outside~$H'$ has an
unknotted subarc outside~$H$, which gets replaced with an isotopic
segment in~$\d H$.  So the arc of~$K^*$ outside~$H'$ is unknotted, as
desired.

Finally, we must consider the case where one endpoint of the arc is
within $H$ and one is outside $H$.  This implies that the half-spaces
are not nested, but divide space into four quadrants,
with $K$ winding once around their intersection line.
The knot $K$ is a connect sum of four arcs, in these quadrants.
The sum of the two arcs outside~$H$ is trivial, as is the sum of the two
outside~$H'$.  So in fact all the knotting happens in the one quadrant
$H\cap H'$.  But now, the arc of~$K^*$ outside~$H'$
consists of the original unknotted arc in one quadrant,
plus a straight segment within $\d H$; thus it is unknotted.
\end{proof}

For the proof of our theorem, we will need to decompose a knot
into its prime pieces, but only when this can be done with a flat plane.

\begin{definition}
We say that a knot $K$ in $\R^3$ is \emph{geometrically composite} if
there is a plane $P$ (cutting $K$ twice) which decomposes $K$ into two
nontrivial connect summands.  Otherwise, we say $K$ is
\emph{geometrically prime}.  Note that a (topologically) prime knot is
necessarily geometrically prime.
\end{definition}

\begin{maintheorem}
If $K$ is any space curve with a knotted inscribed polygon,
then its second hull $\h_2(K)$ is nonempty.
\end{maintheorem}

\begin{proof}
We can replace $K$ by the knotted inscribed polygon,
and by Lemma~\ref{lem:cut} its second hull can only shrink.
If $K$ is now geometrically composite, we can replace it by one of its
two summands, and by Lemma~\ref{lem:cut} again its second hull
will only shrink.  Because the polygonal $K$ is tame,
we need only do this a finite number of times,
and so we may assume that $K$ is geometrically prime.

Let $A:=\bigcap_H \overline{H}$ be the intersection of the
closures of all the half-spaces $H$ essentially containing $K$.
We claim that any plane $P$ through a point $p$ in $A$ cuts $K$ at
least four times.  If not, $P$ would cut $K$ twice.
Because $K$ is geometrically prime, $P$ cannot split
$K$ into two nontrivial summands. So $P$ must be the boundary of an
essential half-space $H$.  We will construct a parallel half-space $H'$
strictly contained in $H$ which still essentially contains $K$. Since
$p \notin \overline{H'}$ and $p \in A$, this is a contradiction.

If $K$ is not contained in~$H$, there are two edges of~$K$ which start
in~$H$ but end outside $H$ (possibly on~$P$).  Let $v$ and $w$ be
their endpoints in~$H$, and $\Gamma$ be
the subarc of $K$ between $v$ and $w$ in~$H$.
If $K \subset H$, let $\Gamma = K$. In either case,
$\Gamma$ is a compact subset of $H$, and so is contained in a smaller
parallel half-space $H'$ which essentially contains $K$. This completes
the proof of the claim.

Since planes through~$A$ cut~$K$ four times,
$A \subset \h_2(K) \subset \h_1(K)$, and we can write
\begin{equation*}
A=\bigcap_H \,\big(\overline H \cap h_1(K)\big).
\end{equation*}
This exhibits $A$ as an infinite intersection of compact sets.  To
show it is nonempty, it suffices to prove that any finite intersection
is nonempty (using the finite intersection property for compact sets).
So suppose $H_1$, \ldots, $H_N$ are half-spaces essentially
containing~$K$.  Clipping the knot successively to these
half-spaces, we define $K_0=K$ and $K_i=\clip_{H_i}(K_{i-1})$.  Since
$K$ is essentially contained in every $H_j$, this is also true
inductively for each of the $K_i$, using Lemma~\ref{lem:clip}.  
Thus $K_N$, a knot isotopic to $K$, is contained in
$$\bigcap_1^N \, \overline{H_i} \cap \h_1(K),$$
showing this to be nonempty, as desired.
\end{proof}

Note that our sets $\overline H\cap \h_1(K)$ are convex.  Therefore, we could
appeal to Helly's theorem instead of the finite intersection property:
as soon as any four of these compact convex sets have nonempty intersection,
they all do.  But the proof given above is no simpler
for the special case $N=4$.

We expect that the set $A$ used in the proof is always a connected
component of $\h_2(K)$.
This essential piece $A\subset\h_2(K)$ is the intersection of
all halfspaces essentially containing~$K$, just as
$h_1(K)$ is the intersection of all halfspaces containing $K$.
We note also that the only place we used the assumption about polygons
was to show that bounding planes of $A$ must cut the knot four times,
and thus that all of $A$, and not merely its (possibly empty) interior,
is in $\h_2(K)$.

Using Milnor's definition of total curvature for arbitrary curves,
our theorem applies to any knot of finite total curvature:

\begin{corollary}
A knotted space curve of finite total curvature has nonempty second hull.
\end{corollary}
\begin{proof}
A lemma of Milnor~\cite{milnor}
shows that any curve of finite total curvature
has an isotopic inscribed polygon, so our main theorem applies.
\end{proof}

\section{Links}
\begin{figure}[ht]
  \begin{center} \incgrh{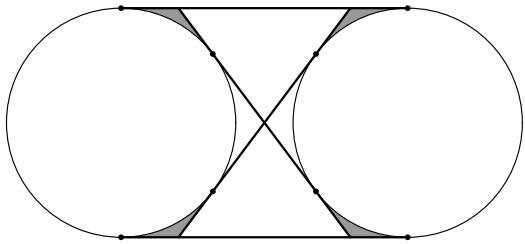}{1.5in} \end{center} 
\caption[The first hull of an unlink.]{
  This figure shows two round circles forming an unlink~$L$ in the plane.
  The convex hull of~$L$ is bounded by a stadium curve.
  Its first hull $\h_1(L)$, however, consists only of the hulls of the two
  components along with the four shaded triangular regions.
  These are bounded by the four lines tangent to both circles.
}\label{fig:link}
\end{figure}

If we apply our definition of $n^{\text{th}}$ hull to links $L$, that
is, to unions of closed curves in space, then it is no longer
necessarily true that $\h_1(L)$ is the entire convex hull of~$L$.
However, Lemma~\ref{lem:cvx} shows that $\h_1(L)$ contains the convex
hull of each component.  In Figure~\ref{fig:link} we show an
example where both inclusions are strict.

For links we can easily find points in higher hulls by intersecting
the hulls of different components.  The next lemma is immediate from
the definitions:

\begin{lemma}\label{lem:union}
If $A$ and $B$ are links, then
$$\h_m(A) \intersect \h_n(B) \subset \h_{m+n}(A\union B)$$
for all $m,n>0$.
\qed
\end{lemma}

This allows us to extend our theorem to the case of links.

\begin{corollary}
Any nontrivial link (of finite total curvature) has a nonempty second hull.
\end{corollary}
\begin{proof}
If the link contains a knotted component, we use the theorem.
Otherwise it must contain two components $A$ and $B$ whose convex
hulls intersect, and we use the last lemma.
\end{proof}

For a Hopf link whose two components are plane curves, Lemma~\ref{lem:union}
is sharp: the second hull consists exactly of the segment along which the
convex hulls of the components intersect. This shows that the second hull
of a link can have zero volume. We conjecture that this never happens
for knots.

\section{Cone Angles}
We originally became interested in second hulls while attempting
to improve lower bounds on the ropelength of knots~\cite{CKS2}.
We wanted to find points in space from which a knot has large cone angle,
or visual angle.

\begin{lemma}\label{lem:cone}
If $p\in \h_n(K)$ then the cone from $p$ to $K$ has cone angle at
least $2\pi n$.
\end{lemma}
\begin{proof}
Consider the radial projection $K'$ of $K$ to the unit sphere around
$p$.  By the definition of $\h_n$, $K'$ intersects each great circle
at least $2n$ times.  But the length of a spherical curve, using
integral geometry as in the proof of Fenchel's theorem, equals $\pi$
times its average number of intersections with great circles~(see
\cite{milnor}).  Thus the length of $K'$, which is the cone angle of
$K$ at $p$, is at least $2\pi n$.
\end{proof}

\begin{corollary}
If $K$ is a nontrivial knot or link (of finite total curvature), there is
some point $p$ from which $K$ has cone angle at least $4\pi$.  \qed
\end{corollary}

For our application to ropelength, however, we needed to know more.
If $K$ has thickness $\tau$, that is, if
a normal tube of radius $\tau$ around $K$ is embedded, then the
cone point~$p$ can be chosen outside this tube.  We will see below, using
quadrisecants, that for generic prime knots, the second hull does extend
outside this thick tube, as desired.

However, a different argument gives $4\pi$ cone points for arbitrary links:
Brian White has a version of the monotonicity theorem for minimal surfaces
with boundary~\cite{white-hild}.  As he points out in work with
Ekholm and Wienholtz~\cite[Thm.~1.3]{EWW},
it shows directly that if we span $K$ with a minimal surface,
and $p$ is a point where this surface intersects itself,
then the cone angle from $p$ to $K$ is at least $4\pi$.
In fact, this result was originally proved in 1983 by
Gromov~\cite[Thm.~8.2.A]{GromovFRM}.

The set of points from which the cone angle to a curve $K$ is at least $2\pi$
has been called the \emph{visual hull}, and is a superset of the convex hull.
(See~\cite{GromovFRM,Gage}.)
We might call the set where cone angle is at least $4\pi$ the
\emph{second visual hull}; as we have mentioned,
this set contains our second hull.
Gromov's result, reproved by our corollary above, is that
any knot has a nonempty second visual hull.

\begin{figure}[ht]
  \begin{center} 
	\vbox{\parbox{2.5in}{\incgrh{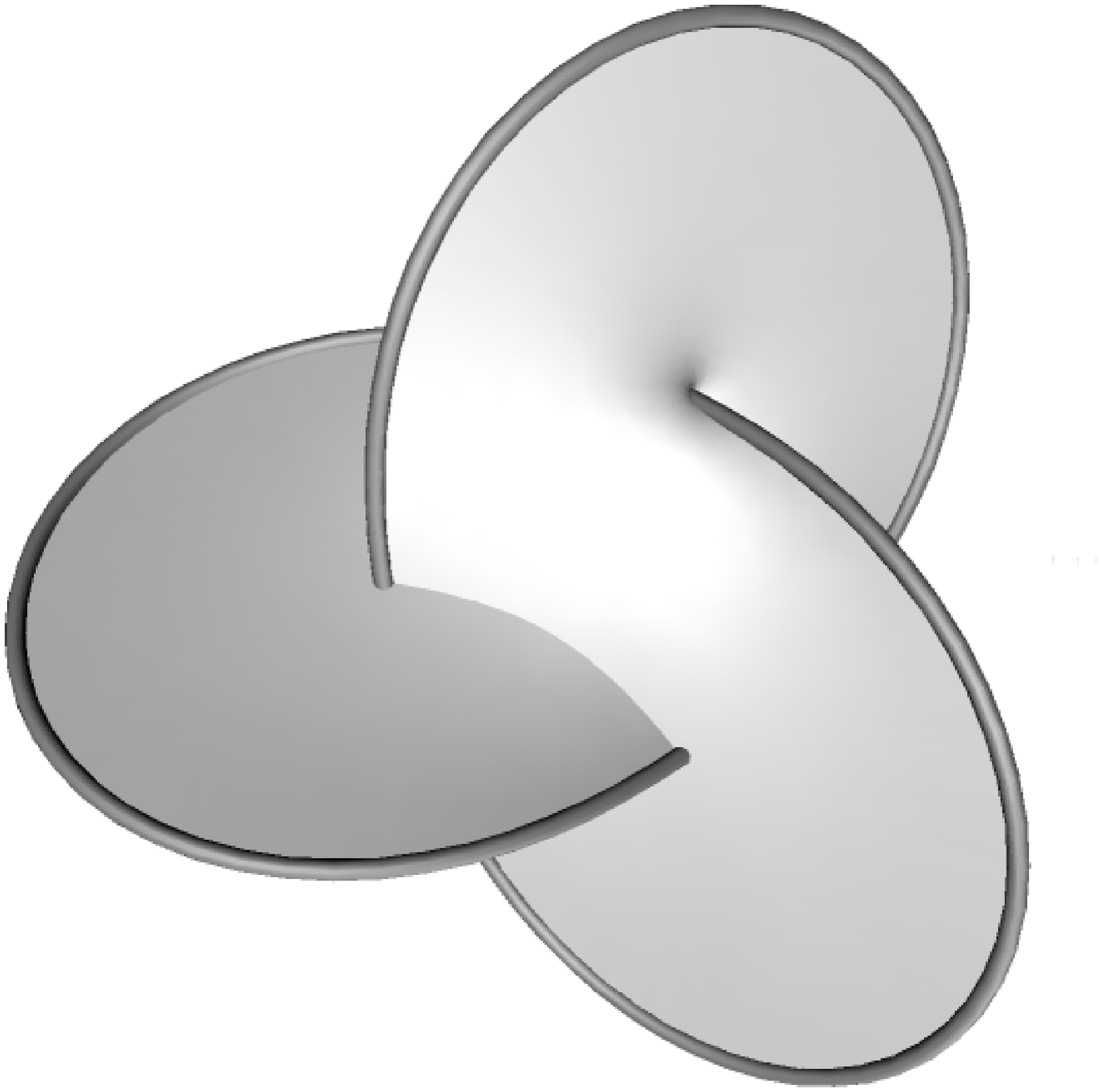}{2in}}}
	\vspace{.2in}
	\vbox{\parbox{2.6in}{\incgrh{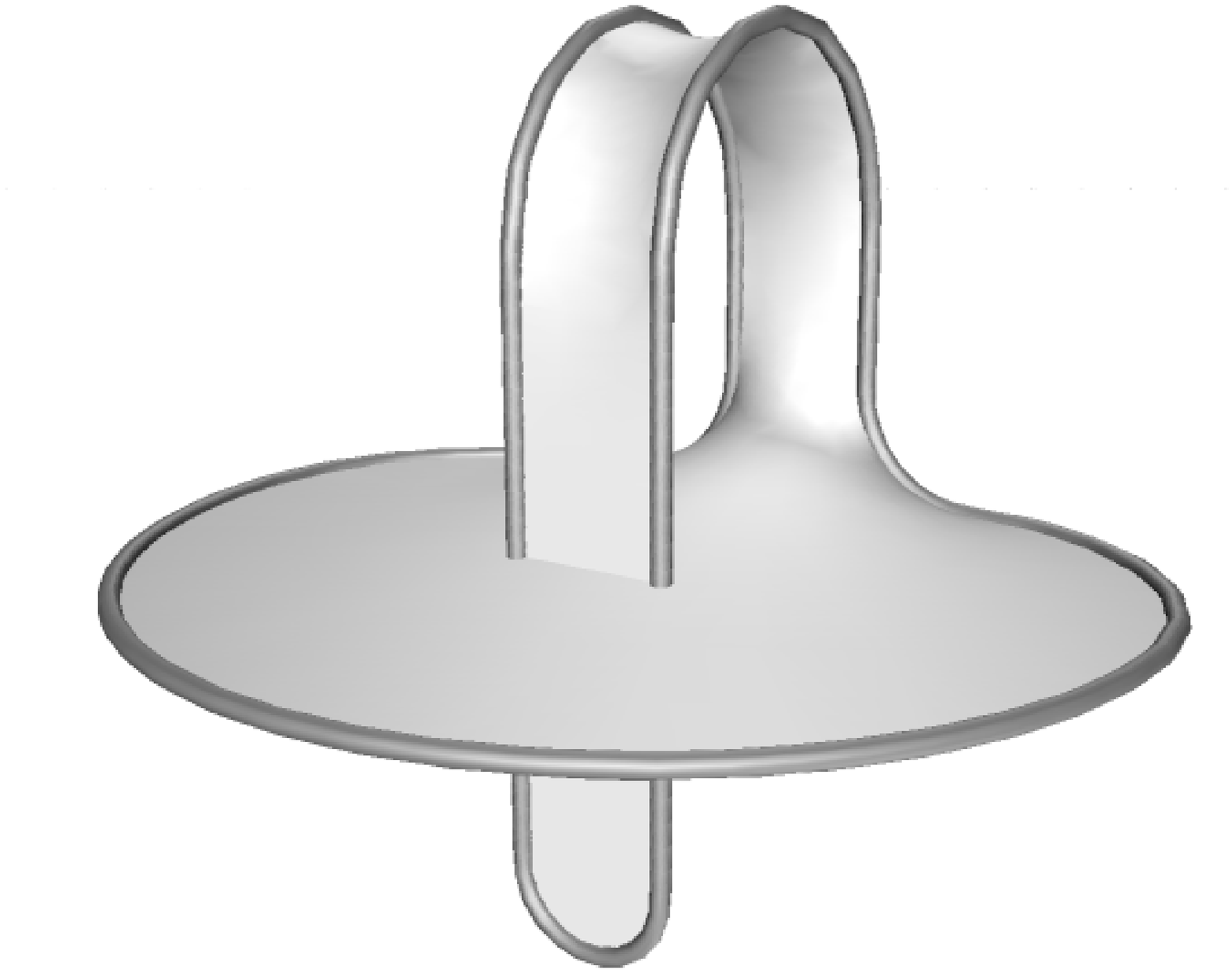}{138bp}}}
   	\end{center} 	
	
	\caption[Constructing $4\pi$ cone points with minimal surfaces.]
	{This figure shows two immersed minimal disks, one
	bounded by a trefoil knot and the other bounded by an
	unknot. By White's version of the monotonicity theorem for
	minimal surfaces, the cone angle to the boundary
        from any point on the line
	of self-intersection of either disk is at least $4\pi$. In the
	trefoil above, this line is contained within the second
	hull. For the unknotted curve below, the second hull is
	empty, and so the line of self-intersections is in the second
        visual hull, but not in the second hull.\label{fig:minimal}}
\end{figure}

We have shown~\cite{CKS2} that if $K$ is a thick knot, then the
self-intersection curve of a minimal spanning disk must go outside the
tube of radius $\tau$, as shown in Figure~\ref{fig:minimal},
giving us what we needed: a point far from the knot with large cone angle.
We note that self-intersection points of a
minimal surface need not be in the second hull of its boundary curve,
as demonstrated by the second example in Figure~\ref{fig:minimal},
though they are necessarily in the second visual hull.

\section{Bridge Number and Total Curvature}
\label{sec:hullbridgetotal}

\begin{figure}[ht]
  \begin{center} \includegraphics[height=1.8in]{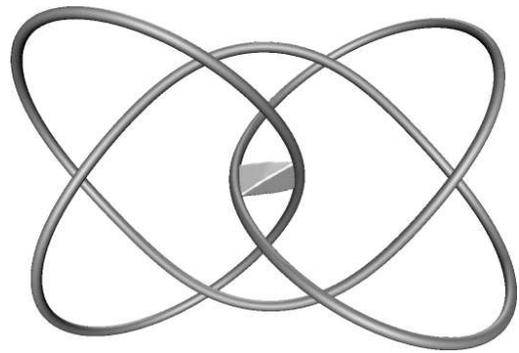} \end{center}
	\caption[The third hull of $8_{19}$.]
	{This presumed M\"obius-energy-minimizing configuration
	(from~\cite{KS-knot}) of the $(3,4)$-torus knot $8_{19}$,
	has a small but nonempty third hull.
	We do not know whether every configuration of $8_{19}$
	has nonempty third hull.
	This knot and the button knot $8_{18}$ are the first
	prime $3$-bridge knots in Rolfsen's table~\cite{rol}.
	As in Figure~\ref{fig:trefoil2hull}, we show a numerical
	approximation of the third hull, and we expect that the true
	hull is somewhat smoother.  \label{fig:8_19}}
\end{figure}

We can define the \emph{hull number} of a link to be the largest~$n$
for which the $n^{\text{th}}$ hull is nonempty (and the hull number of
a link type to be the minimum over all representatives).  Perhaps
sufficiently complicated link types have hull number greater than two,
but we know of no way to prove this.  We can get an easy upper bound:
\begin{proposition}
The hull number of any link is bounded above by its bridge number.
\end{proposition}
\begin{proof}
To say $L$ is a link of bridge number $n$ means there is
some height function with just $n$ minima and $n$ maxima.
Thus $L$ will have empty $(n+1)^{\text{st}}$ hull,
because it cuts no horizontal plane more than $2n$ times.
\end{proof}
Note that we don't always have equality, at least for composite links.
For instance, the connect sum of two 2-bridge links (like Hopf links
or trefoil knots) has bridge number three.  But obvious
configurations, where each summand is in its 2-bridge presentation and
they are well separated vertically, have empty third hull.  Perhaps
any prime 3-bridge knot has nonempty third hull as in Figure~\ref{fig:8_19}:
certainly there are
planes in every direction which are cut six times by the knot.
To extend our proof to this case, we might need to consider
the notion of thin position for knots,
as introduced by Gabai~\cite{Gabai}
and studied by Thompson~\cite{Thompson-thin}.

Milnor's version~\cite{milnor} of the F{\'a}ry/Milnor theorem says
that the total curvature of a knot is greater than $2\pi$ times its
bridge number; by our proposition, the same is true replacing bridge
number by hull number. Since our main theorem says hull
number is at least two for a nontrivial knot, this provides an
alternate proof of F{\'a}ry/Milnor.

The following lemma, which also appears in~\cite{EWW},
when combined with Lemma~\ref{lem:cone}, gives yet another
proof of F{\'a}ry/Milnor.
\begin{lemma}
The cone angle of a knot $K$, from any point, is at most the total
curvature of $K$.
\end{lemma}
\begin{proof}
Applying Gauss-Bonnet to the cone itself, we see that the cone angle
equals the total geodesic curvature of $K$ in the cone, which is no
greater than the total curvature of $K$ in space.
\end{proof}

Kuiper~\cite{Kui-sb} defined the \emph{superbridge number} $\sb(K)$
of a curve $K$ as half the maximum number of times $K$ cuts any plane.
Like our notion of \nthhull/, this is invariant under affine
transformations of the ambient space.  Kuiper observed that
the total curvature of $K$ is bounded above by $2\pi\sb(K)$,
and that there is a sequence of affine transformations of $K$ (stretched
in the direction perpendicular to any plane cut $2\sb(K)$ times)
whose total curvature approaches that limit $2\pi\sb(K)$.
We note that, viewed from any point near the center of mass, these
stretched curves also have cone angles approaching the same value $2\pi\sb(K)$.

F{\'a}ry's proof~\cite{fary} of the F{\'a}ry/Milnor theorem
showed that the total curvature of a space curve is
the average total curvature of its planar projections,
and that any knot projection has total curvature at least $4\pi$.
In fact, his proof of the latter statement actually showed
that any knot projection has nonempty second hull.

When we first started to think about second hulls, we rederived
this result, not knowing F{\'a}ry's proof.  We hoped that perhaps
a curve whose projections all had nonempty second hulls would
necessarily also have one, but never succeeded in this line of argument.
The next lemma does give a sort of converse, relating hulls and projections.
\begin{lemma}
If $K$ is a closed curve in $\R^d$, and $\Pi$ is any orthogonal
projection to $\R^{d-1}$, then $\Pi(\h_n(K))\subset \h_n(\Pi(K))$.
\end{lemma}
\begin{proof}
Call the projection direction ``vertical''.  If $P$ is a hyperplane in
$\R^{d-1}$ through $p\in\Pi(\h_n(K))$, then the vertical hyperplane
$\Pi^{-1}P$ in $\R^d$ passes through a point in $\h_n(K)$, so $K$ cuts
it $2n$ times.  Thus $\Pi(K)$ cuts $P$ at least $2n$ times.
\end{proof}

This lemma explains why, for instance, the projected third hull
in Figure~\ref{fig:8_19} must lie entirely within the central
bigon of the knot projection.

\section{Quadrisecants and Second Hulls}

Any knot or link in space has a quadrisecant,
that is, a straight line which intersects the link four times.
This was first proved by Pannwitz~\cite{Pan} for generic polygonal links;
the full result is due to Kuperberg~\cite{Kup}.
Quadrisecants can be used to prove a thickened version of our theorem.
In particular, the middle segment of a quadrisecant for $K$
often lies in the second hull of $K$.

In the case of links, if there are two components $A$ and $B$ with
nonzero linking number, Pannwitz~\cite{Pan} (and later Morton and
Mond~\cite{MorMon}) showed that there is an $ABAB$ quadrisecant,
meaning one where the components are seen in that order along the
secant line.  The entire mid-segment then lies in the intersection of
the convex hulls of $A$ and $B$, and thus in $\h_2(A\cup B)$.  For the
thickened version, simply note that this mid-segment starts inside the
tube around $A$ and ends inside the disjoint one around $B$, so it
must leave the tubes altogether in the middle.

For knots, we can compare the linear ordering of the four
intersection points (along the quadrisecant line) with their circular
ordering (along the knot).
Viewed from a point on the mid-segment of a quadrisecant,
the four intersection points lie (two each) at the north
and south poles of the visual sphere.
If these intersection points alternate between the north and south
poles, we call the quadrisecant \emph{alternating}.  The mid-segment
of an alternating quadrisecant for $K$ will be contained
in the second hull of $K$: if we project $K$ to the unit sphere
around a point on the mid-segment, the projected curve cuts every great
circle four times, since it visits the poles in alternating order.
Thus the following conjecture would give an alternate proof of our main
theorem.

\begin{conjecture}
Any nontrivial knot has an alternating quadrisecant.
\end{conjecture}

Elizabeth Denne has already made good progress towards proving this
conjecture as part of her doctoral dissertation at the University of Illinois.
She makes use of ideas from Pannwitz~\cite{Pan}, Kuperberg~\cite{Kup}
and Schmitz~\cite{Schmitz}.

For generic polynomial curves in prime knot classes,
we can prove a thickened version of our
theorem using Kuperberg's notion of an essential (or
``topologically nontrivial'') quadrisecant~\cite{Kup}.
\begin{definition}
A secant $S$ of $K$ is \emph{trivial} if its endpoints are on the same
component of $K$, dividing it into arcs $X$ and $Y$,
and if one of the two circles $S\cup X$, $S\cup Y$
spans a disk whose interior avoids $K$.  The disk
may intersect $S$ and itself.  A quadrisecant is \emph{essential} if neither
the mid-segment nor either end-segment is trivial. 
\end{definition}

\begin{figure}[ht]
  \begin{center} \includegraphics[height=144bp]{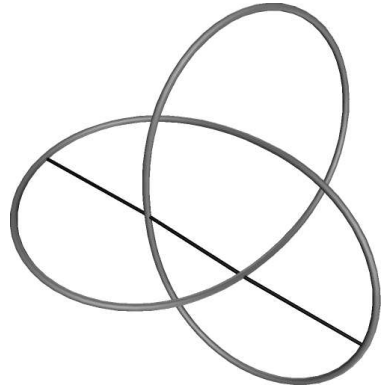} \end{center}
    \caption[A quadrisecant of a trefoil knot.]
    {This figure shows an essential quadrisecant of a
    particular configuration of a trefoil knot. By Lemma~\ref{lem:quad},
    the mid-segment of this quadrisecant is contained within
    the second hull of~$K$. In fact, this quadrisecant is alternating
    as well, so our earlier argument gives an alternate proof that
    the mid-segment is in the second hull of $K$. Since this is the
    same configuration of trefoil knot used in
    Figure~\ref{fig:trefoil2hull}, it is interesting to compare
    this picture with the image of the entire second hull.}
    \label{fig:quad}
\end{figure}

We start by proving a lemma about secants of unknots.

\begin{lemma}
\label{lem:unknot}
Let $K$ be an unknot. Any secant of $K$ is trivial on either side.
\end{lemma}

\begin{proof}
Suppose that $S$ is a secant of $K$.  Pick either one of the arcs of $K$
joining the endpoints of $S$. We will show that the closed curve $C$ formed
by joining this arc to $S$ bounds a disk disjoint from $K$. We begin
the disk by constructing an embedded ribbon bounded by $C$ and a
parallel curve, $C'$, chosen so that the linking number of $K$ and
$C'$ is zero. By the definition of linking number, this implies that
$C'$ is homologous to zero in the complement of $K$. Since $K$ is the
unknot, $C'$ is therefore \emph{homotopic} to zero in the complement
of~$K$.  Thus, $C'$ bounds a disk disjoint from $K$. Joining this disk
to the ribbon completes the proof.
\end{proof}

\begin{lemma}
\label{lem:quad}
Let $K$ be a geometrically prime knot.  The mid-segment of any
essential quadrisecant for $K$ is contained in the second hull of $K$.
\end{lemma}

\begin{proof}
Let $Q$ be a quadrisecant of $K$; if its mid-segment is not contained
in the second hull, there is a plane $P$ that meets the mid-segment of
$Q$ and meets $K$ only twice.  Because $K$ is geometrically prime, it
is trivial on one side of~$P$.  Let $K^*$ be the unknotted arc on
this side, completed with a straight line segment in $P$.  By
Lemma~\ref{lem:unknot}, the end-segment~$S$ of~$Q$ on this side is
trivial with respect to~$K^*$.  That is, $S$ together with one arc of $K^*$
bounds a disk $D$ whose boundary is disjoint from~$P$.

By design, $D$ does not intersect $K^*$, but it may cross $P$ and
intersect $K\setm K^*$.  We claim that we can simplify $D$ so that it
does not cross $P$. Thus the secant $S$ of $Q$ is inessential with
respect to $K$.

To check the claim, consider~$D$ as a map from the standard disk to~$\R^3$.
In general position, the inverse image of~$P$ is a collection
of circles, and we can cut off the disk at the outermost circles.
We know $D$ is interior-disjoint from all of $K^*$,
including the line segment lying in $P$.
But~$D$ may intersect~$P$ in a complicated loop~$L$ that
encircles the pair of points $K\cap P$ many times.
In this case, after cutting~$D$ along this circle, we attach not a
flat disk in~$P$ but instead a bigger disk that avoids~$K$.
\end{proof}

\begin{lemma}
The mid-segment of an essential quadrisecant for $K$ leaves any embedded
normal tube around $K$.
\end{lemma}
\pagebreak[2]
\begin{proof}
If not, it would be isotopic within the tube
to part of the core curve, that is, to one arc of $K$, contradicting
the definition of essential.
\end{proof}

The result of~\cite{Kup} that a generic polynomial link has an essential
quadrisecant, combined with these lemmas, proves the
thickened version of our theorem, for generic prime knots:

\begin{theorem} \label{thm:quad}
If $K$ is a generic polynomial representative of a prime knot class,
and $T$ is an open embedded normal tube around $K$,
then $\h_2(K) \setm T$ is nonempty.
\qed
\end{theorem}

\section{Other Ambient Spaces}

We have been assuming that our knots lie in Euclidean space $\R^3$.
However, our notion of second hull is not a metric notion but
rather a projective notion, depending only on incidence relations of planes.

It follows that our definition applies in hyperbolic space $\H^3$,
and that any knotted curve in~$\H^3$ has a nonempty second hull:
Simply embed hyperbolic space, in Klein's projective
model, in a ball in $\R^3$.  Then the hyperbolic second
hull of any curve is the Euclidean second hull of its image in the model.

Recently, two proofs have been given showing that 
the F{\'a}ry/Milnor theorem holds for knots in any Hadamard manifold,
that is, any simply connected manifold of nonpositive curvature.
Alexander and Bishop~\cite{AlxBsh:FM} find a sequence of
inscribed polygons limiting to a quadruply covered segment,
while Schmitz~\cite{Schmitz} comes close to constructing
a quadrisecant.  Arguments like either of these could be used
to give alternate proofs of our theorem in Euclidean space.

We thus suspect that there should be some notion of second hull
for curves in Hadamard manifolds.  The problem is that
there is no natural way to extend our definition,
because there is no analog of a plane in this general context.
If some reasonable definition for second hull can be found so
that Lemma~\ref{lem:cut} remains true, then there should be no
problem proving that knots have nonempty second hulls. 

\acknowledgments
We gratefully acknowledge helpful conversations with many
colleagues, including Colin Adams, Stephanie Alexander, Dick Bishop,
Elizabeth Denne, Herbert Edelsbrunner, Mike Gage and Frank Morgan.
This work has been partially supported by the NSF through grants
DMS-99-02397 (to Cantarella), DMS-00-72342 (to Kuperberg),
DMS-00-76085 (to Kusner) and DMS-00-71520 (to Sullivan).

\bibliography{thick} 
\vspace{-3ex}
\end{document}